\documentclass[10pt,a4paper, german]{scrartcl}

\usepackage[german]{babel}
\usepackage[latin1]{inputenc}
\usepackage{scrtime}
\usepackage{amsmath}
\usepackage{amsopn}
\usepackage{amsthm}
\usepackage{amscd}
\usepackage{scrtime}
\usepackage{amsxtra}
\usepackage{upref}
\usepackage{amsfonts}
\usepackage{pifont}
\usepackage{hhline}
\usepackage{color}
\usepackage{amstext}
\usepackage[dvips]{geometry}
\usepackage{graphicx}
\usepackage{epsfig}
\begin{document}

\parbox{1000pt}{
\begin{picture}(0,0)(45,629)
\put(89,-2){\psfig{figure=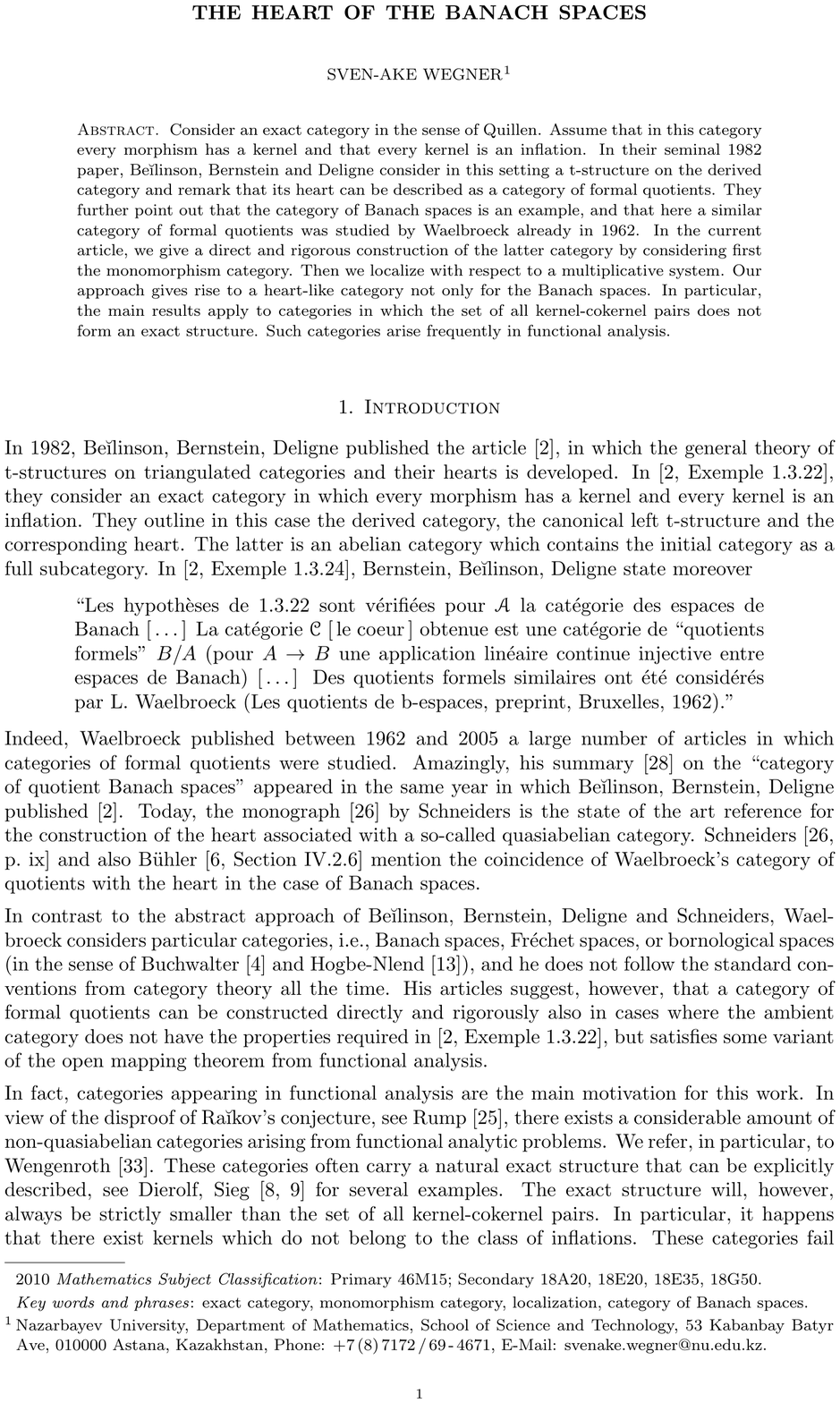, width=21cm}}
\end{picture}
}
\newpage
\parbox{1000pt}{
\begin{picture}(0,0)(45,629)
\put(89,-2){\psfig{figure=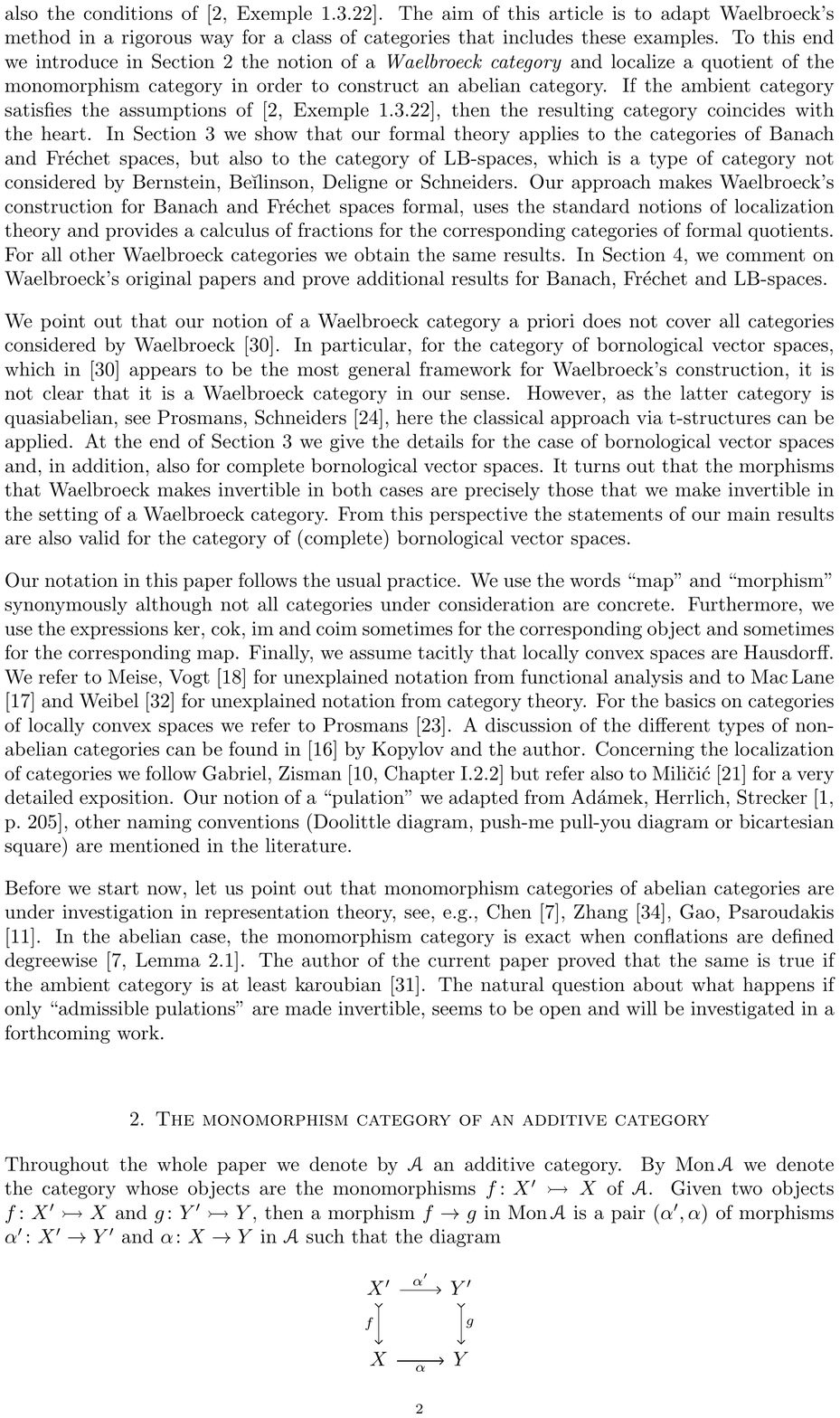, width=21cm}}
\end{picture}
}
\newpage
\parbox{1000pt}{
\begin{picture}(0,0)(45,629)
\put(89,-2){\psfig{figure=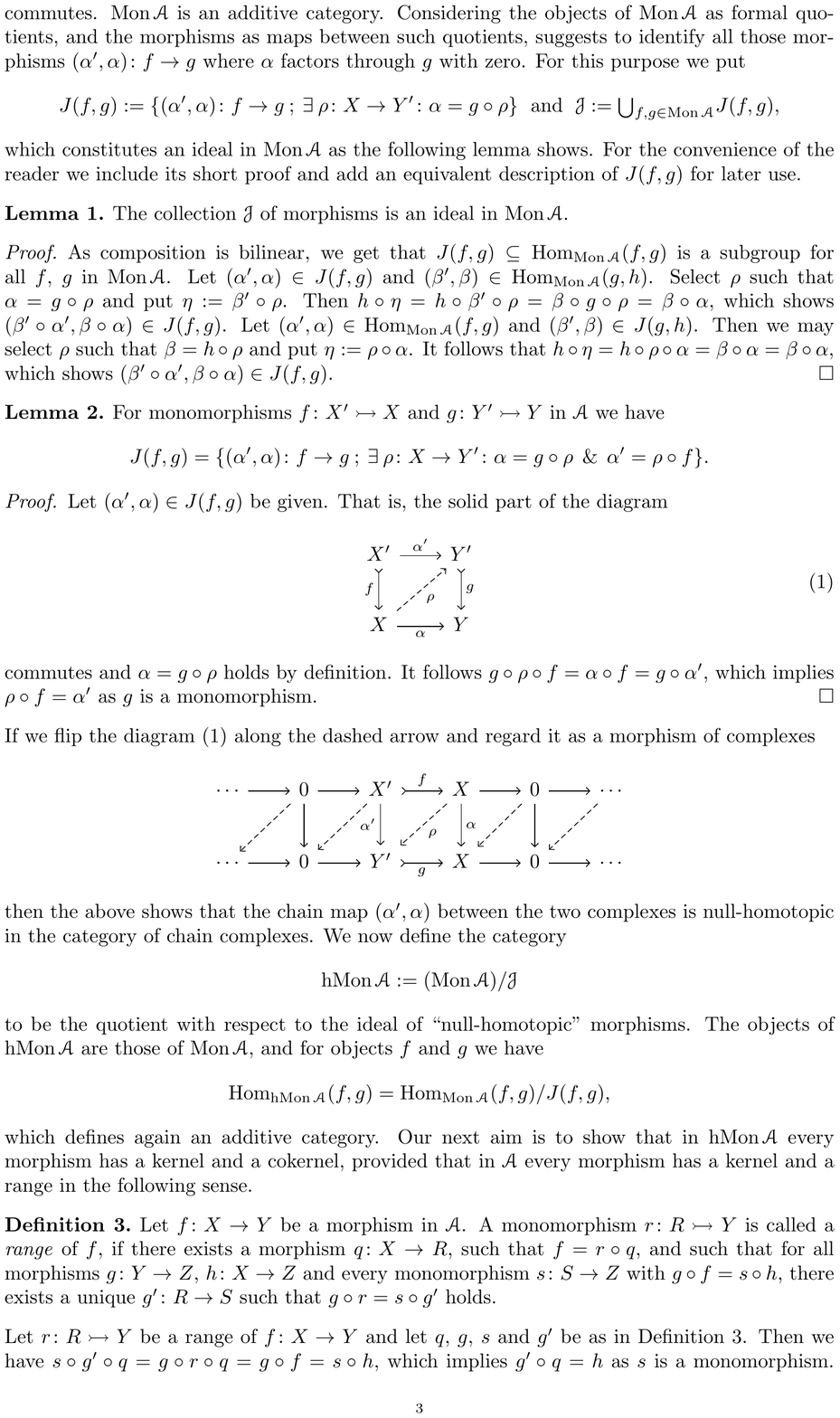, width=21cm}}
\end{picture}
}
\newpage
\parbox{1000pt}{
\begin{picture}(0,0)(45,629)
\put(89,-2){\psfig{figure=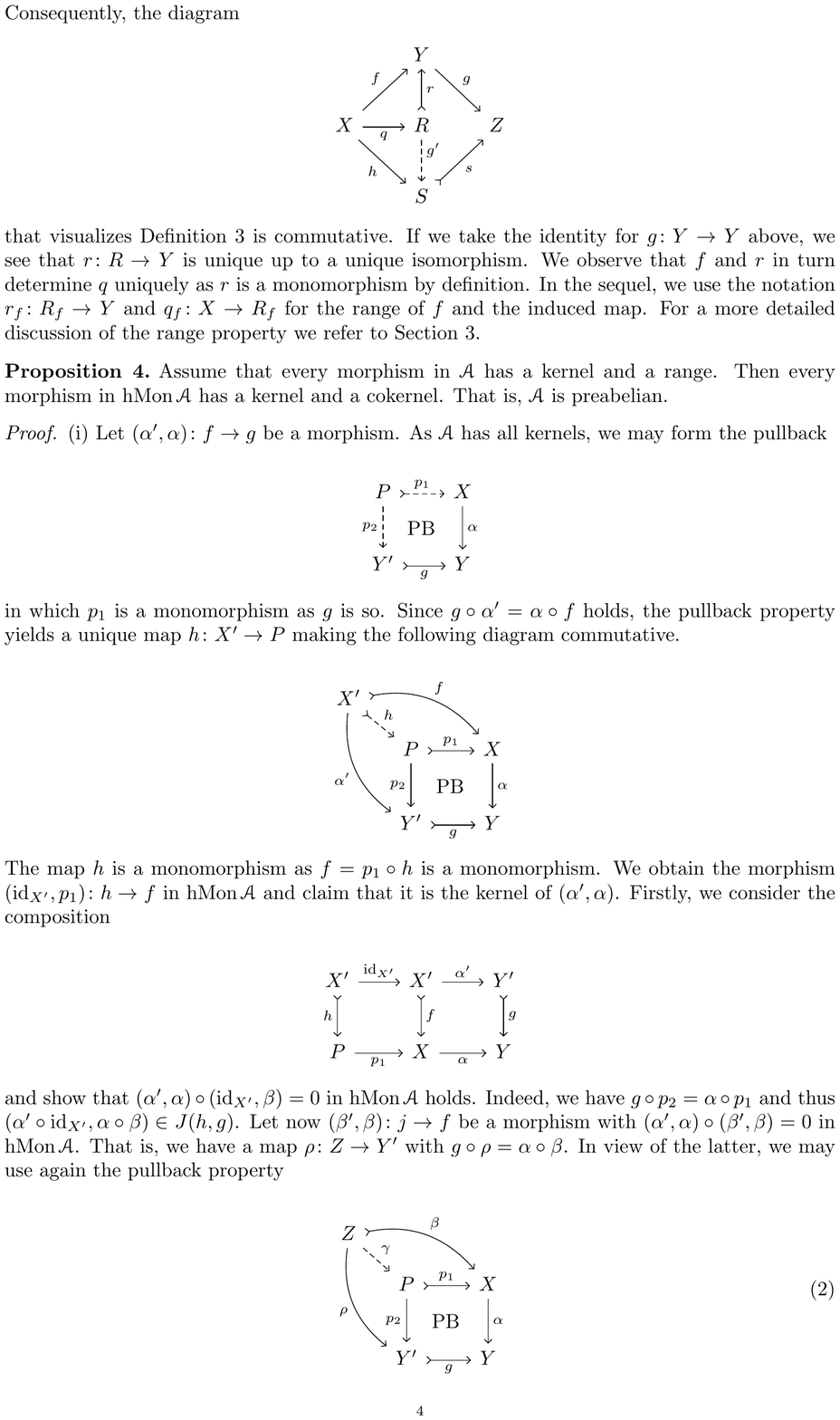, width=21cm}}
\end{picture}
}
\newpage
\parbox{1000pt}{
\begin{picture}(0,0)(45,629)
\put(89,-2){\psfig{figure=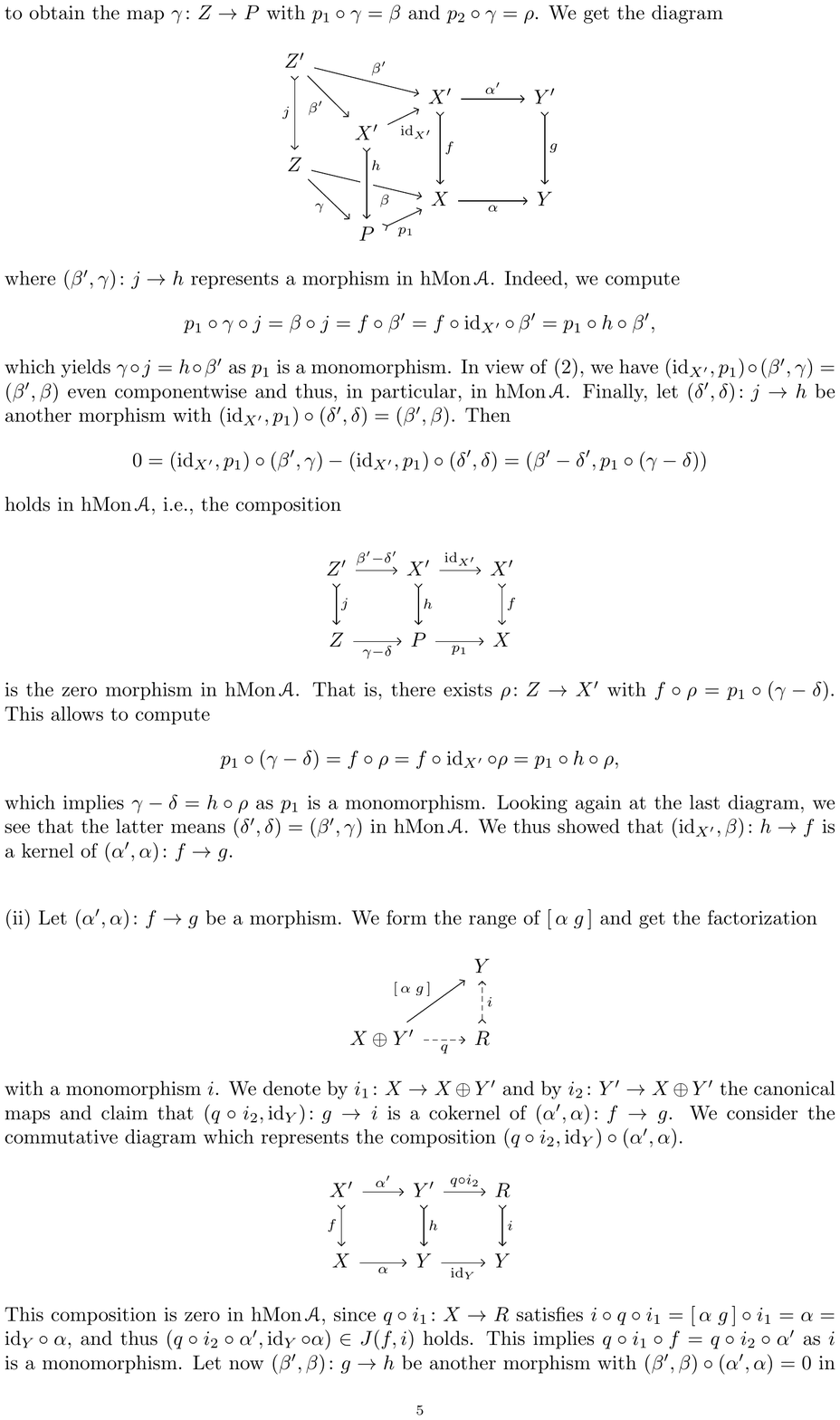, width=21cm}}
\end{picture}
}
\newpage
\parbox{1000pt}{
\begin{picture}(0,0)(45,629)
\put(89,-2){\psfig{figure=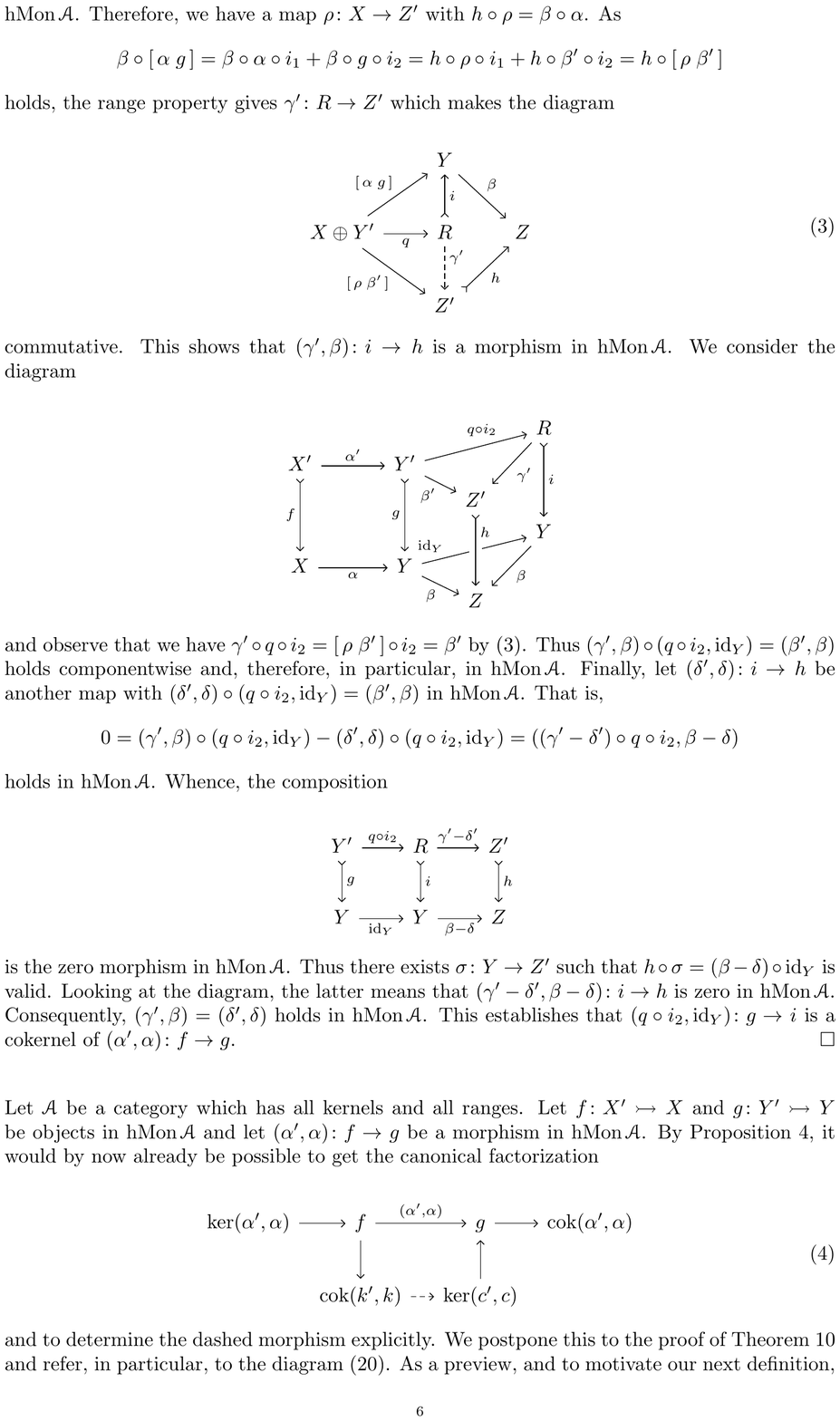, width=21cm}}
\end{picture}
}
\newpage
\parbox{1000pt}{
\begin{picture}(0,0)(45,629)
\put(89,-2){\psfig{figure=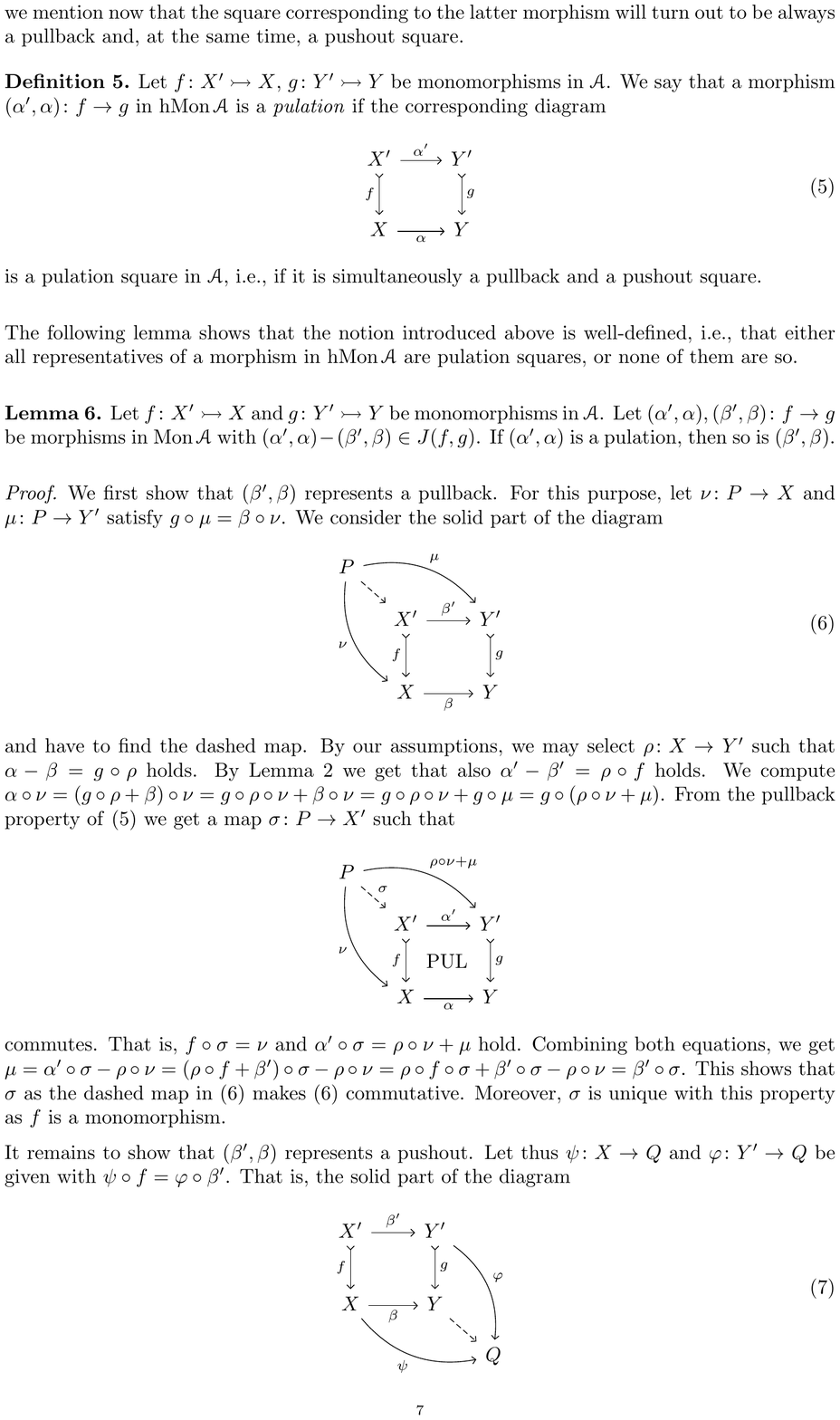, width=21cm}}
\end{picture}
}
\newpage
\parbox{1000pt}{
\begin{picture}(0,0)(45,629)
\put(89,-2){\psfig{figure=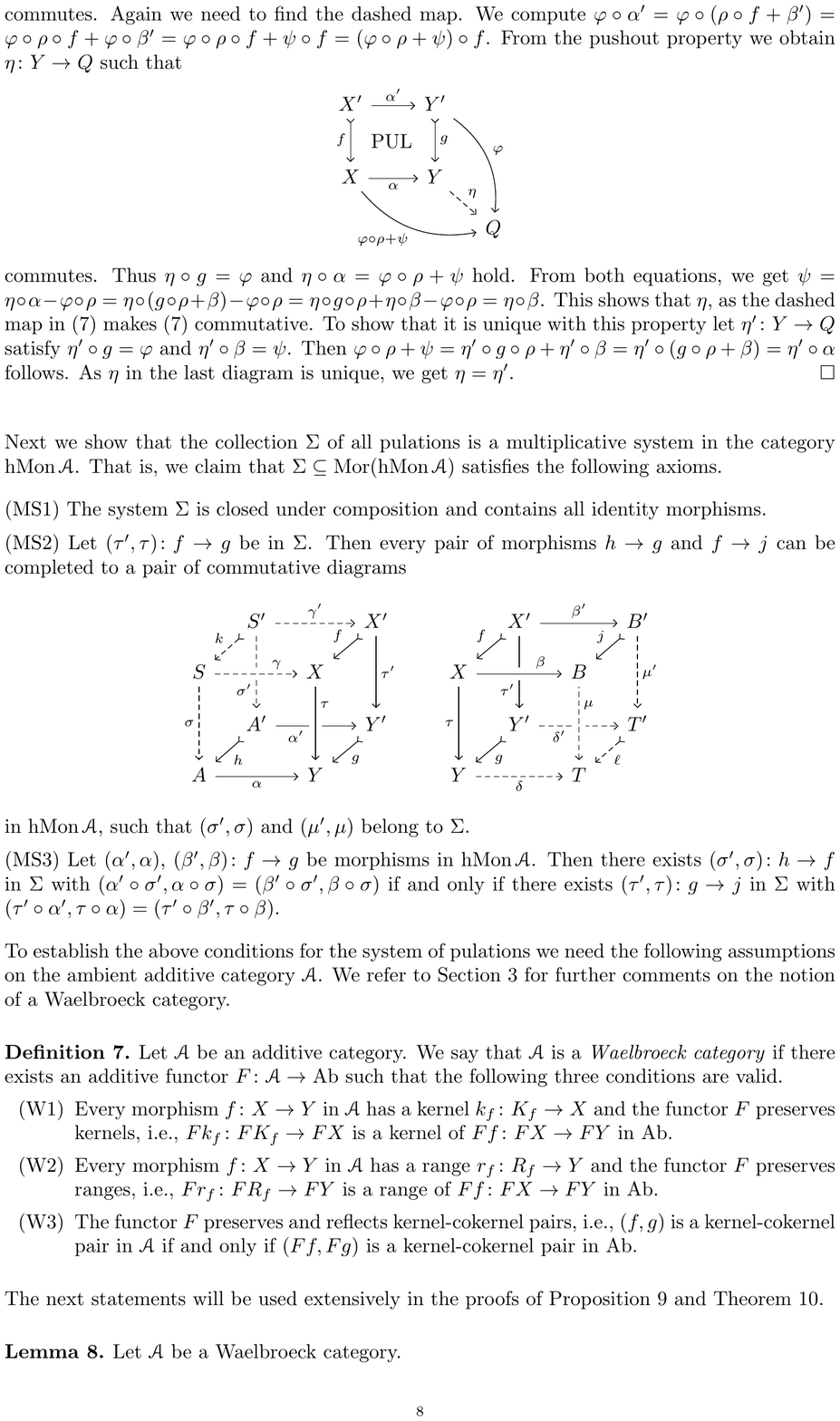, width=21cm}}
\end{picture}
}
\newpage
\parbox{1000pt}{
\begin{picture}(0,0)(45,629)
\put(89,-2){\psfig{figure=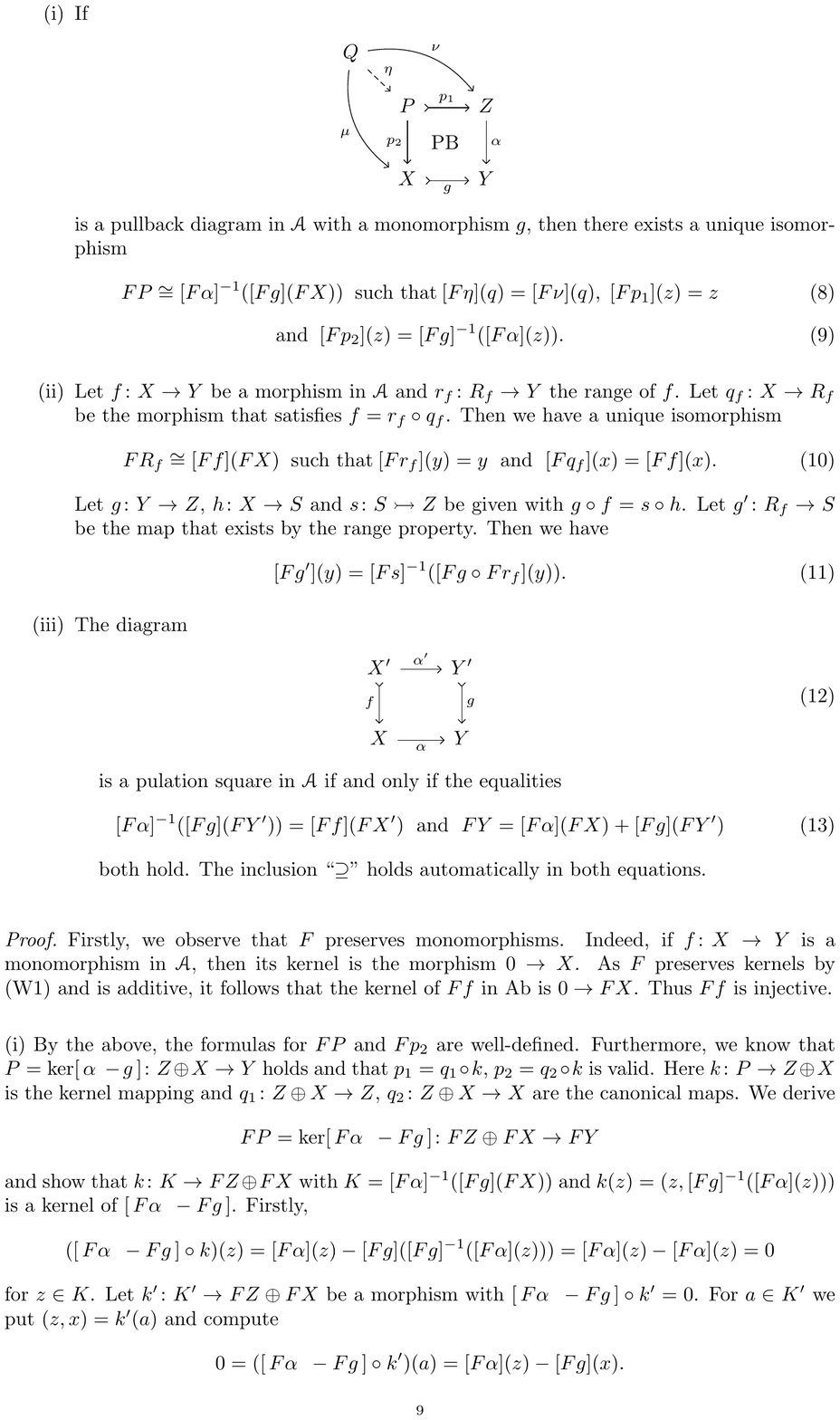, width=21cm}}
\end{picture}
}
\newpage
\parbox{1000pt}{
\begin{picture}(0,0)(45,629)
\put(89,-2){\psfig{figure=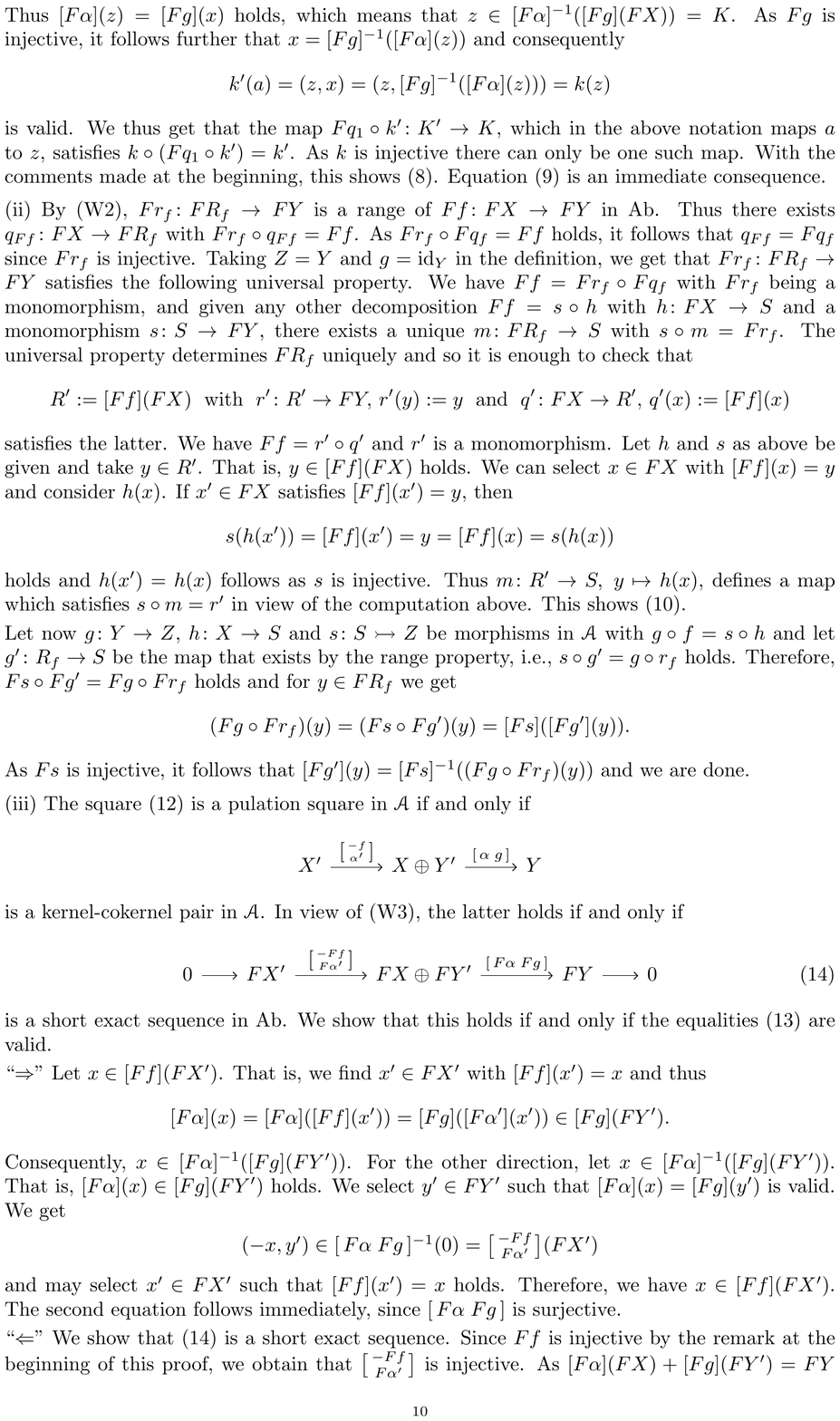, width=21cm}}
\end{picture}
}
\newpage
\parbox{1000pt}{
\begin{picture}(0,0)(45,629)
\put(89,-2){\psfig{figure=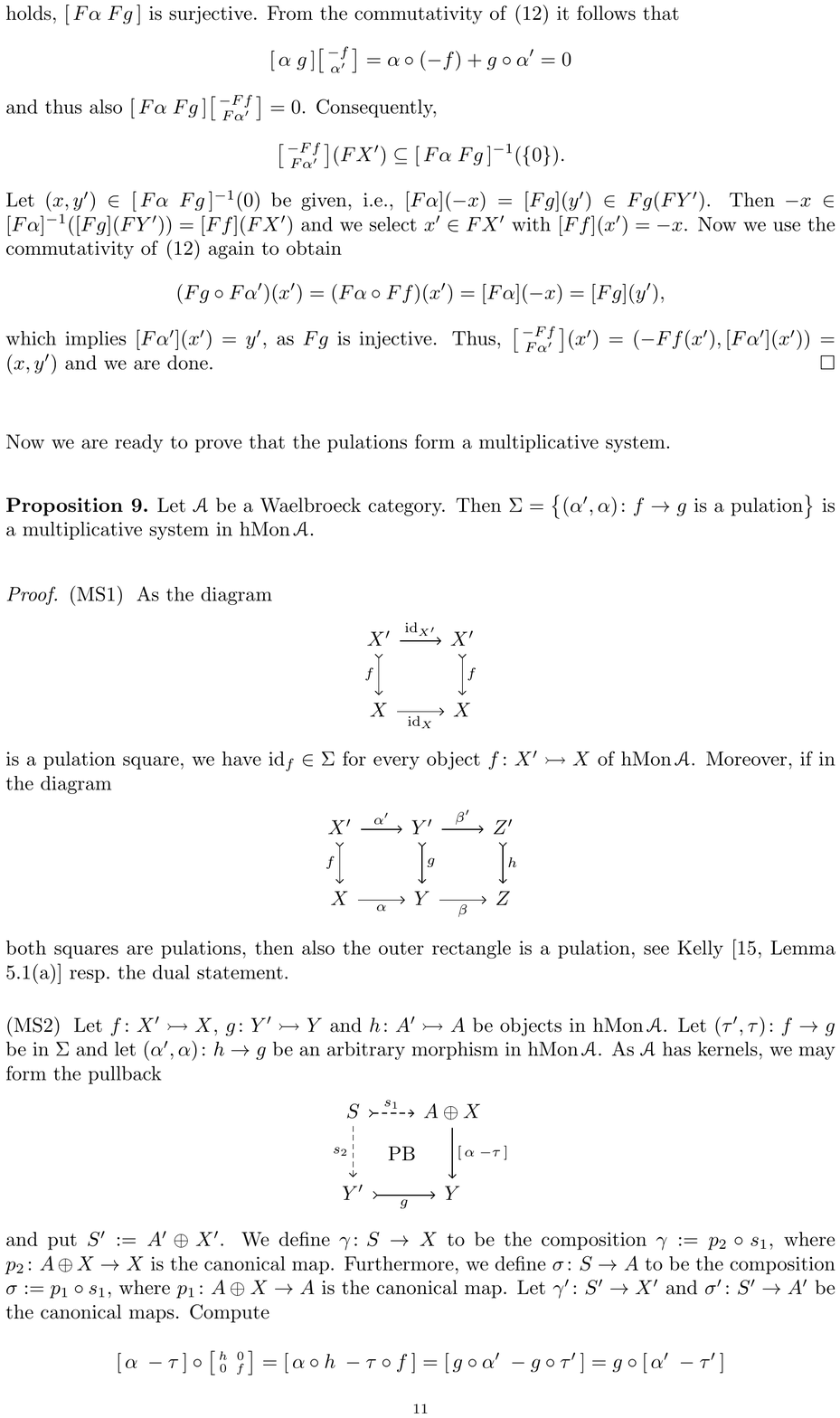, width=21cm}}
\end{picture}
}
\newpage
\parbox{1000pt}{
\begin{picture}(0,0)(45,629)
\put(89,-2){\psfig{figure=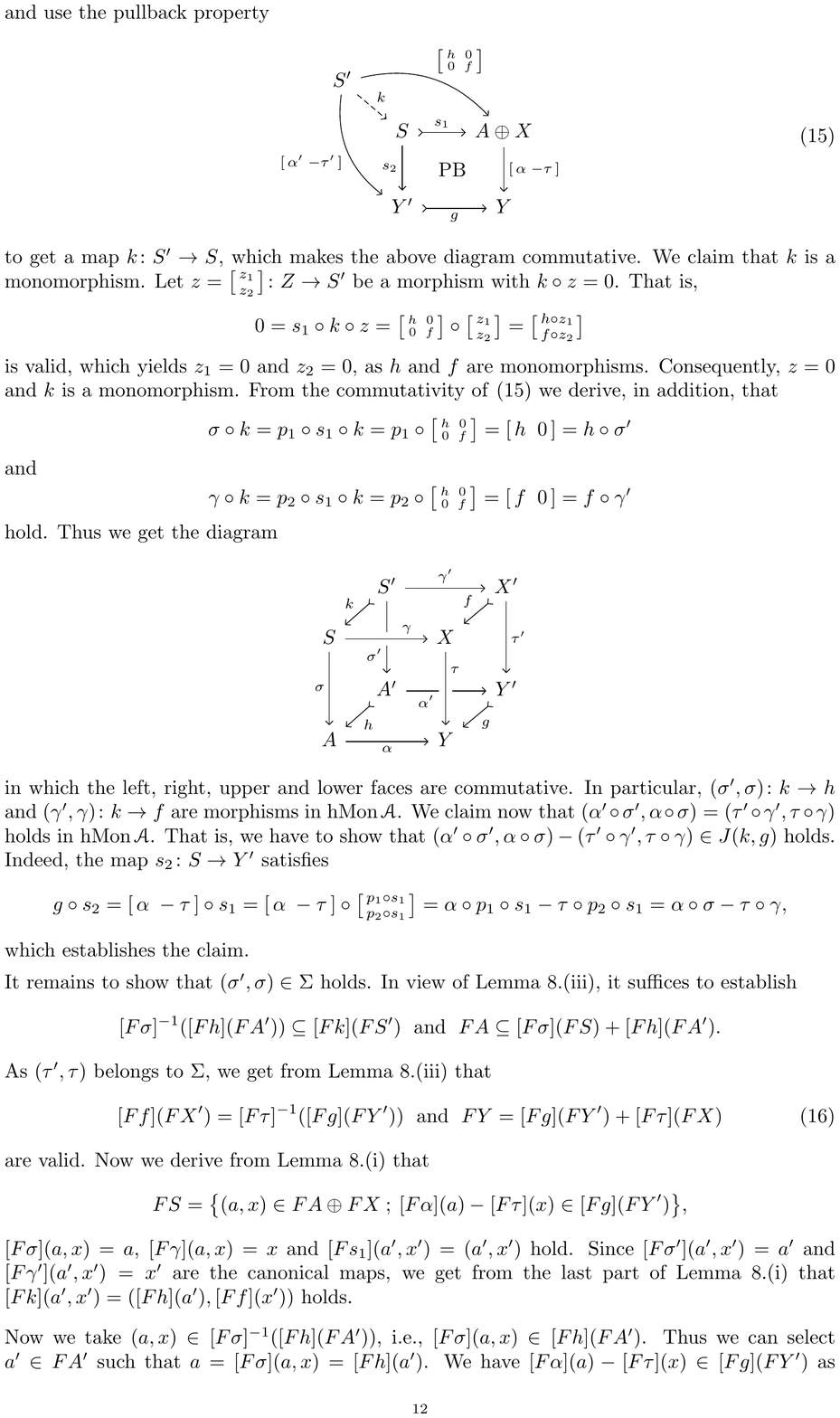, width=21cm}}
\end{picture}
}
\newpage
\parbox{1000pt}{
\begin{picture}(0,0)(45,629)
\put(89,-2){\psfig{figure=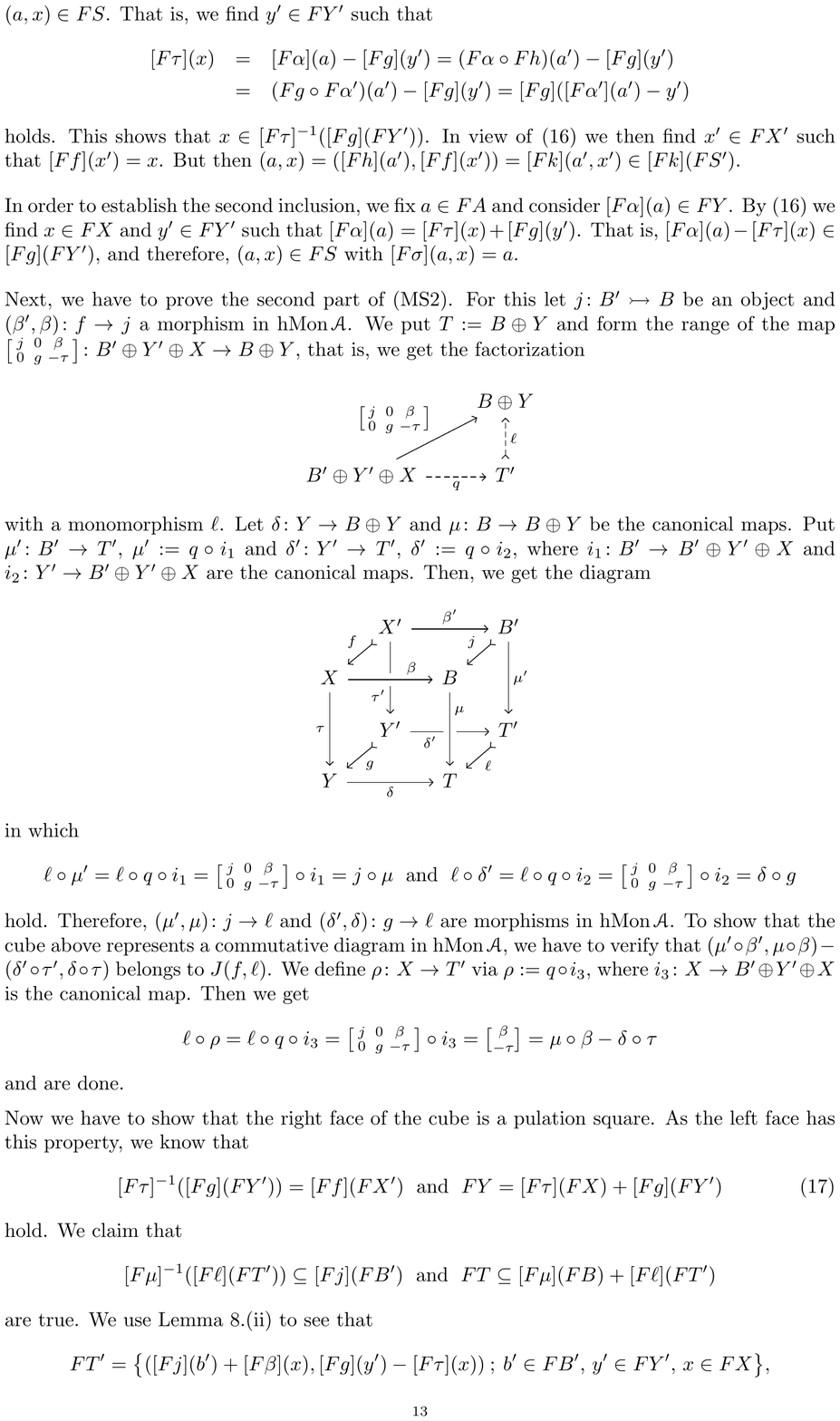, width=21cm}}
\end{picture}
}
\newpage
\parbox{1000pt}{
\begin{picture}(0,0)(45,629)
\put(89,-2){\psfig{figure=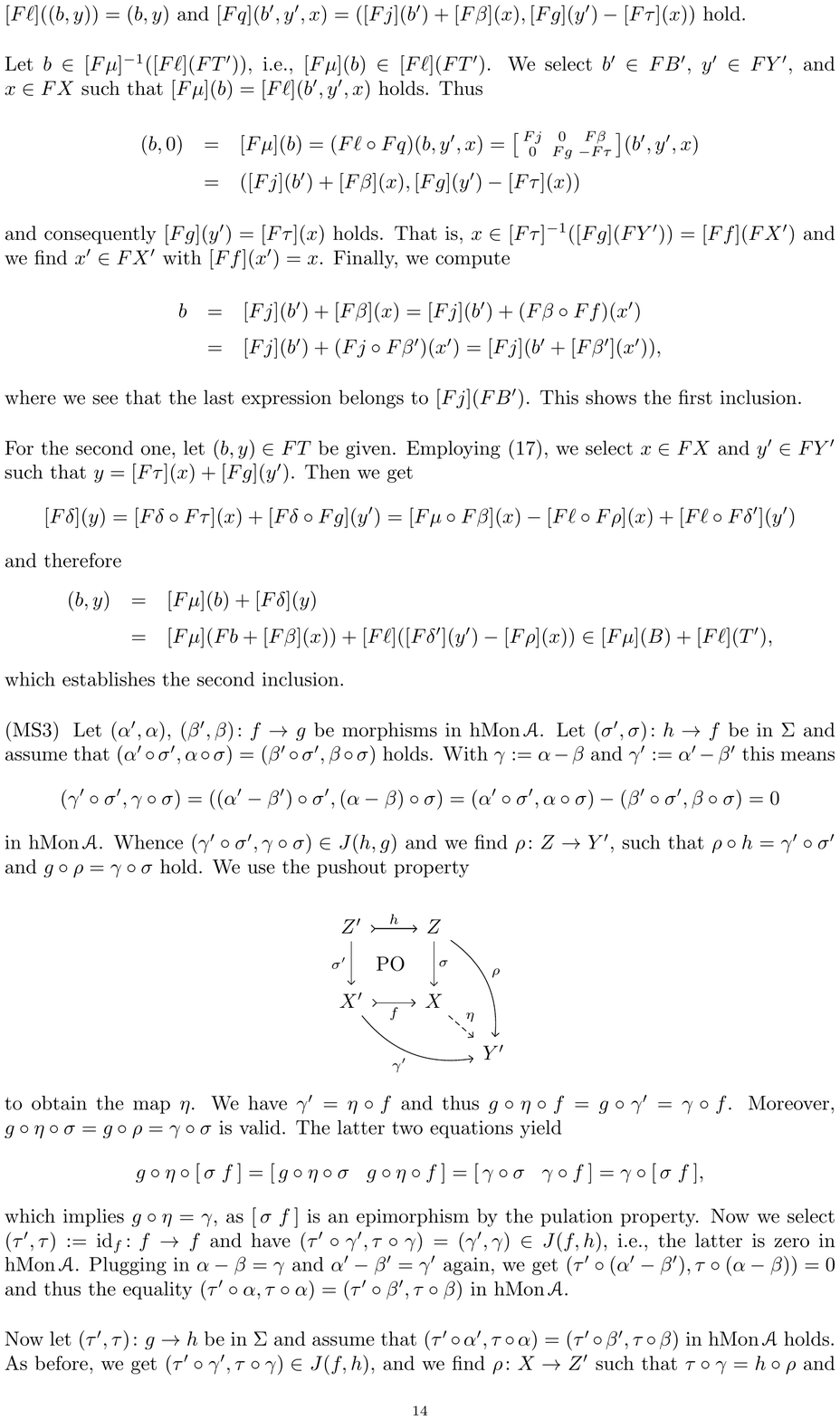, width=21cm}}
\end{picture}
}
\newpage
\parbox{1000pt}{
\begin{picture}(0,0)(45,629)
\put(89,-2){\psfig{figure=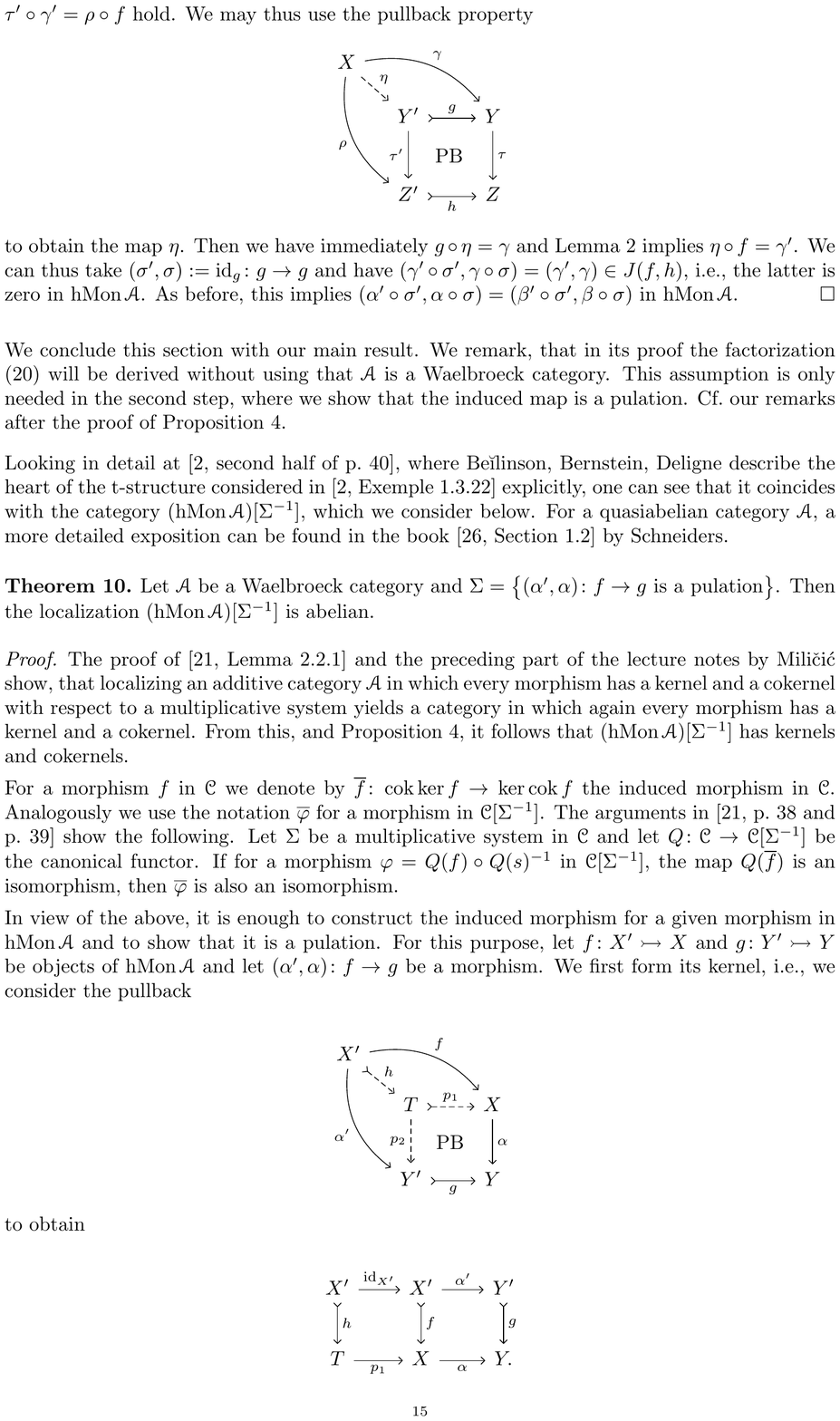, width=21cm}}
\end{picture}
}
\newpage
\parbox{1000pt}{
\begin{picture}(0,0)(45,629)
\put(89,-2){\psfig{figure=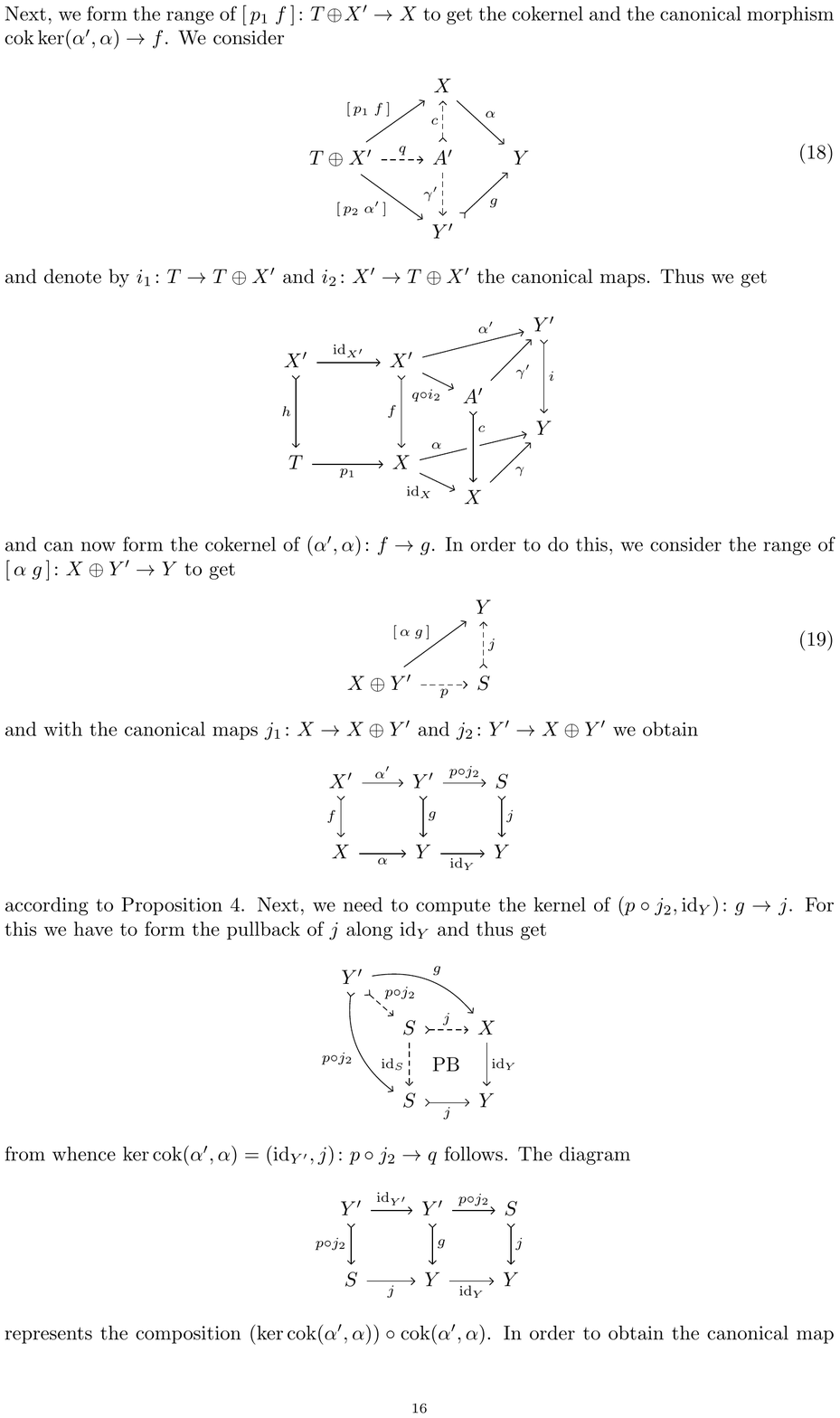, width=21cm}}
\end{picture}
}
\newpage
\parbox{1000pt}{
\begin{picture}(0,0)(45,629)
\put(89,-2){\psfig{figure=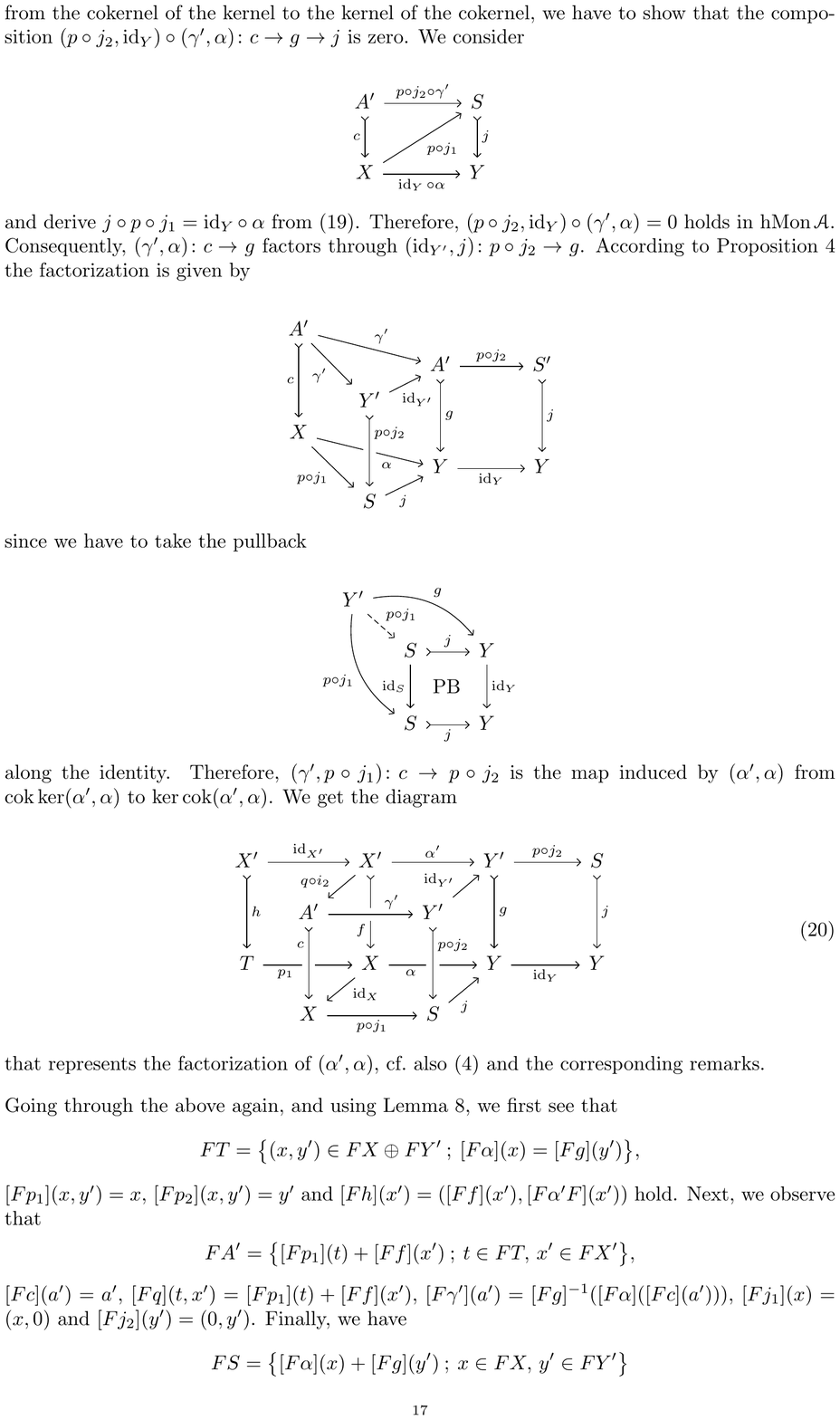, width=21cm}}
\end{picture}
}
\newpage
\parbox{1000pt}{
\begin{picture}(0,0)(45,629)
\put(89,-2){\psfig{figure=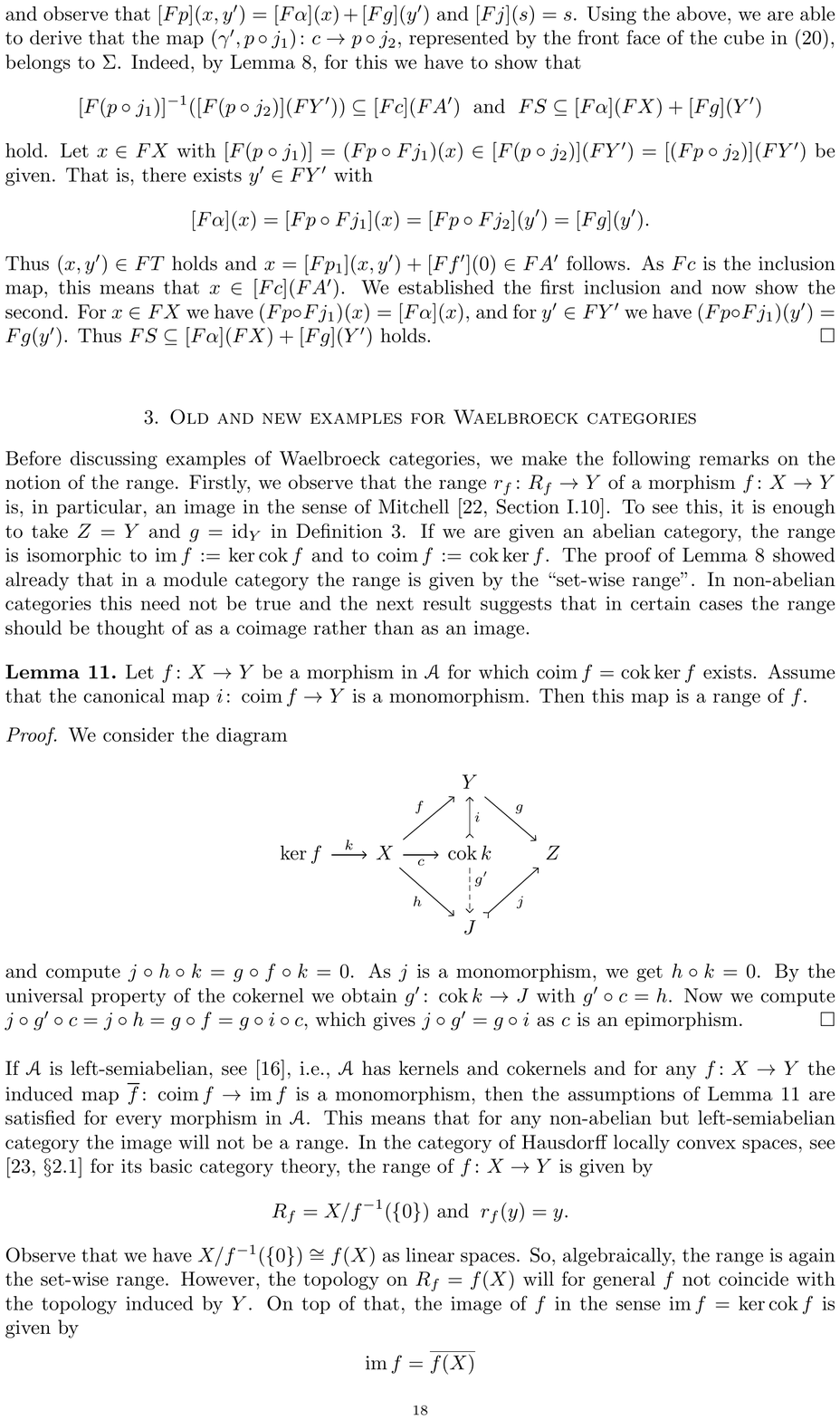, width=21cm}}
\end{picture}
}
\newpage
\parbox{1000pt}{
\begin{picture}(0,0)(45,629)
\put(89,-2){\psfig{figure=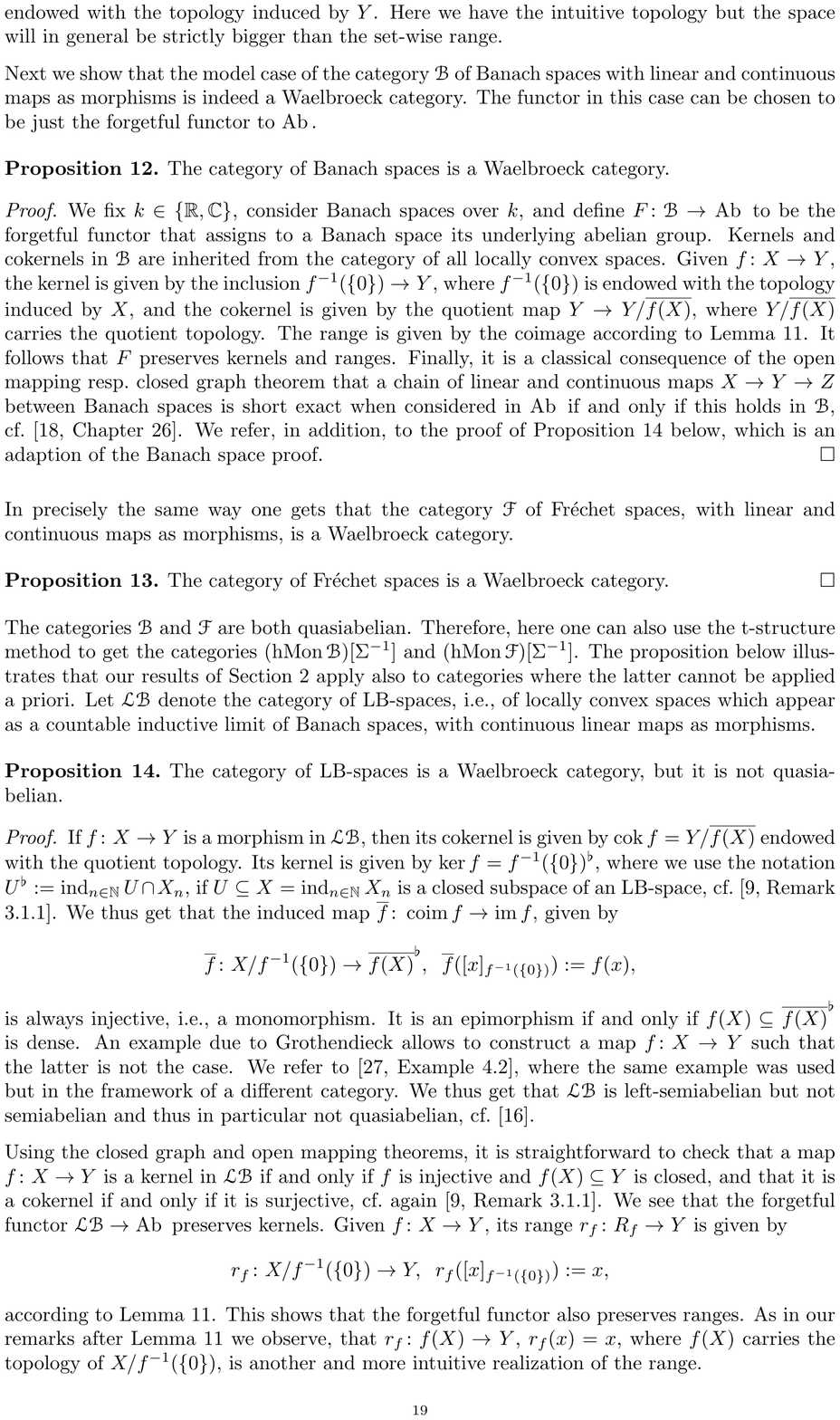, width=21cm}}
\end{picture}
}
\newpage
\parbox{1000pt}{
\begin{picture}(0,0)(45,629)
\put(89,-2){\psfig{figure=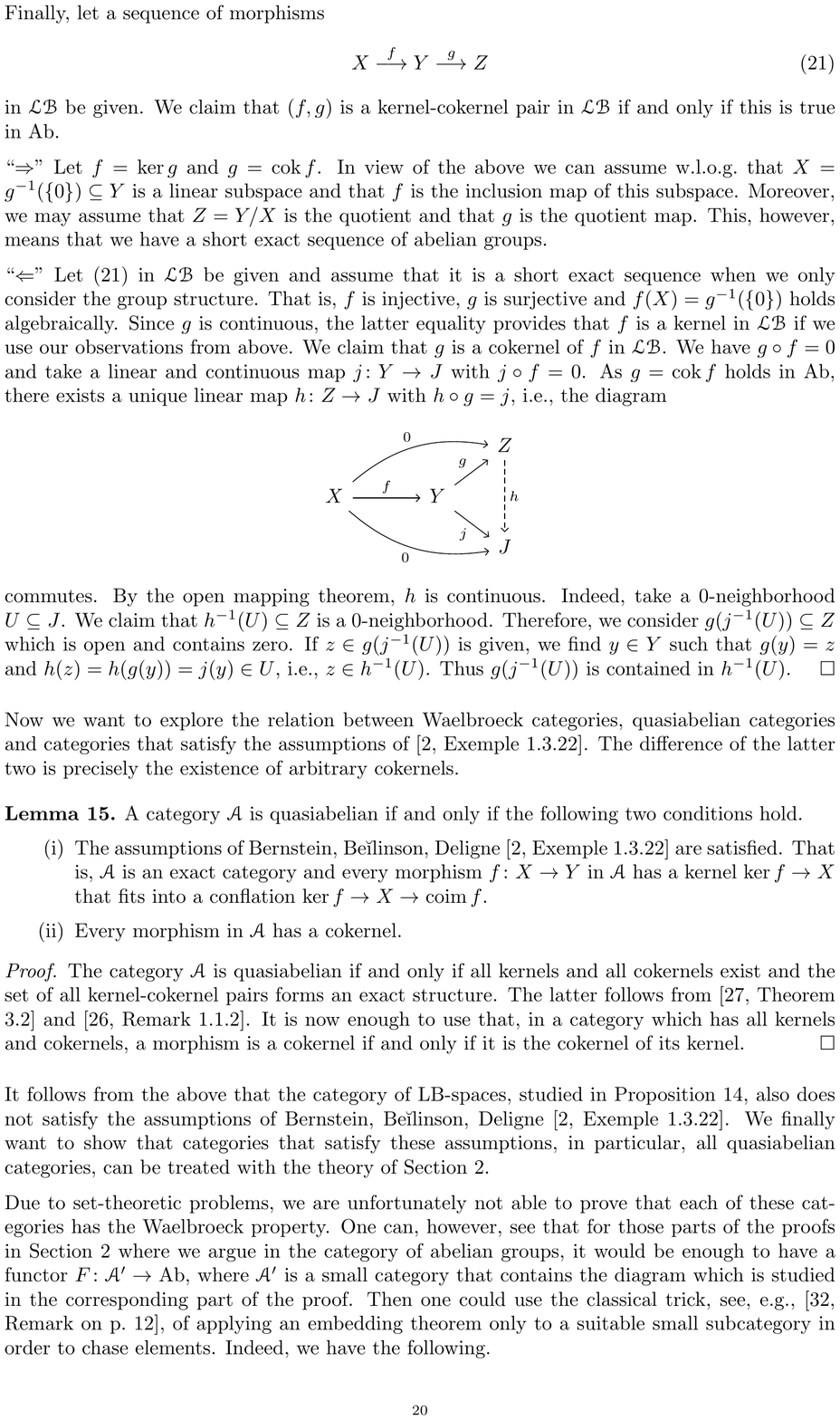, width=21cm}}
\end{picture}
}
\newpage
\parbox{1000pt}{
\begin{picture}(0,0)(45,629)
\put(89,-2){\psfig{figure=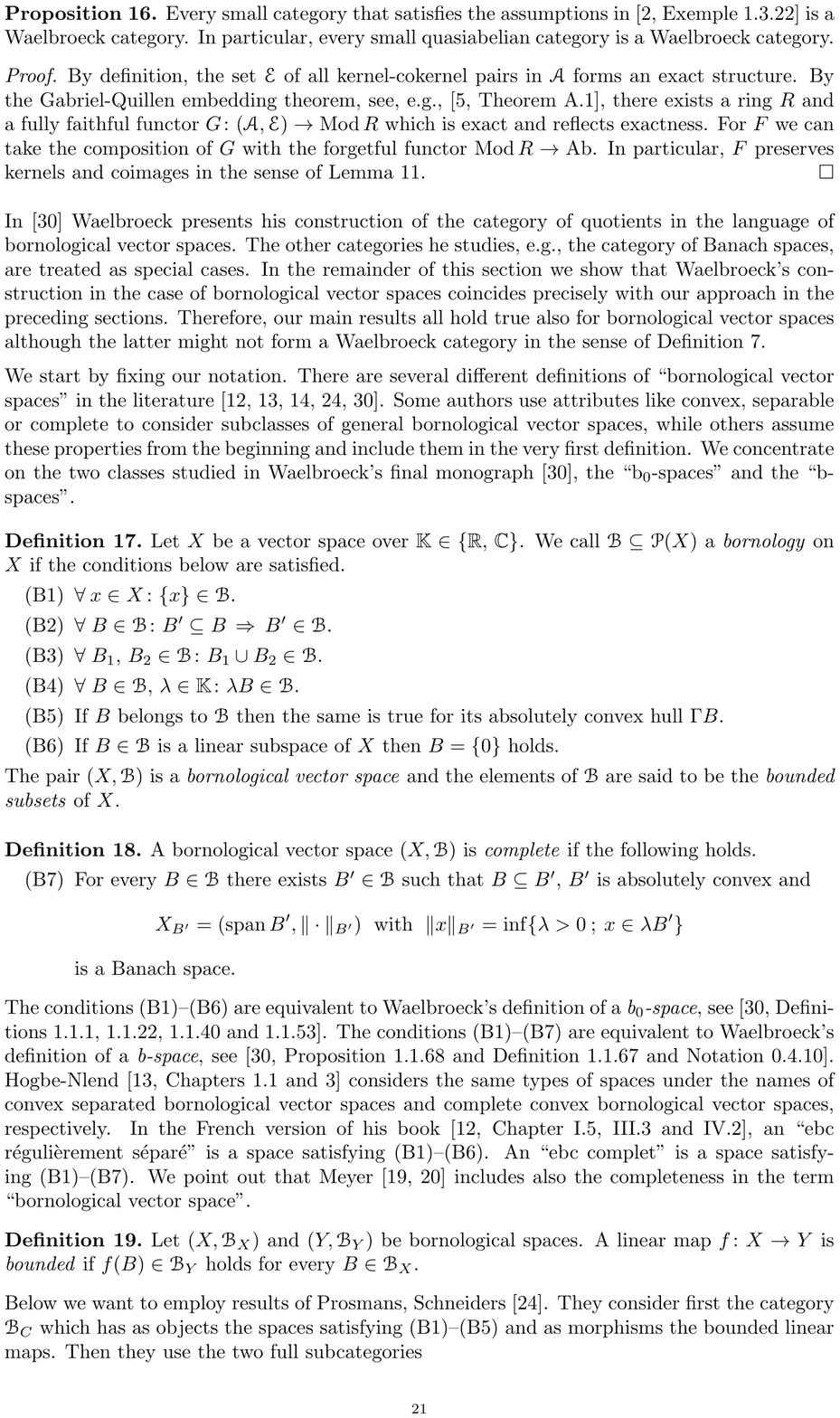, width=21cm}}
\end{picture}
}
\newpage
\parbox{1000pt}{
\begin{picture}(0,0)(45,629)
\put(89,-2){\psfig{figure=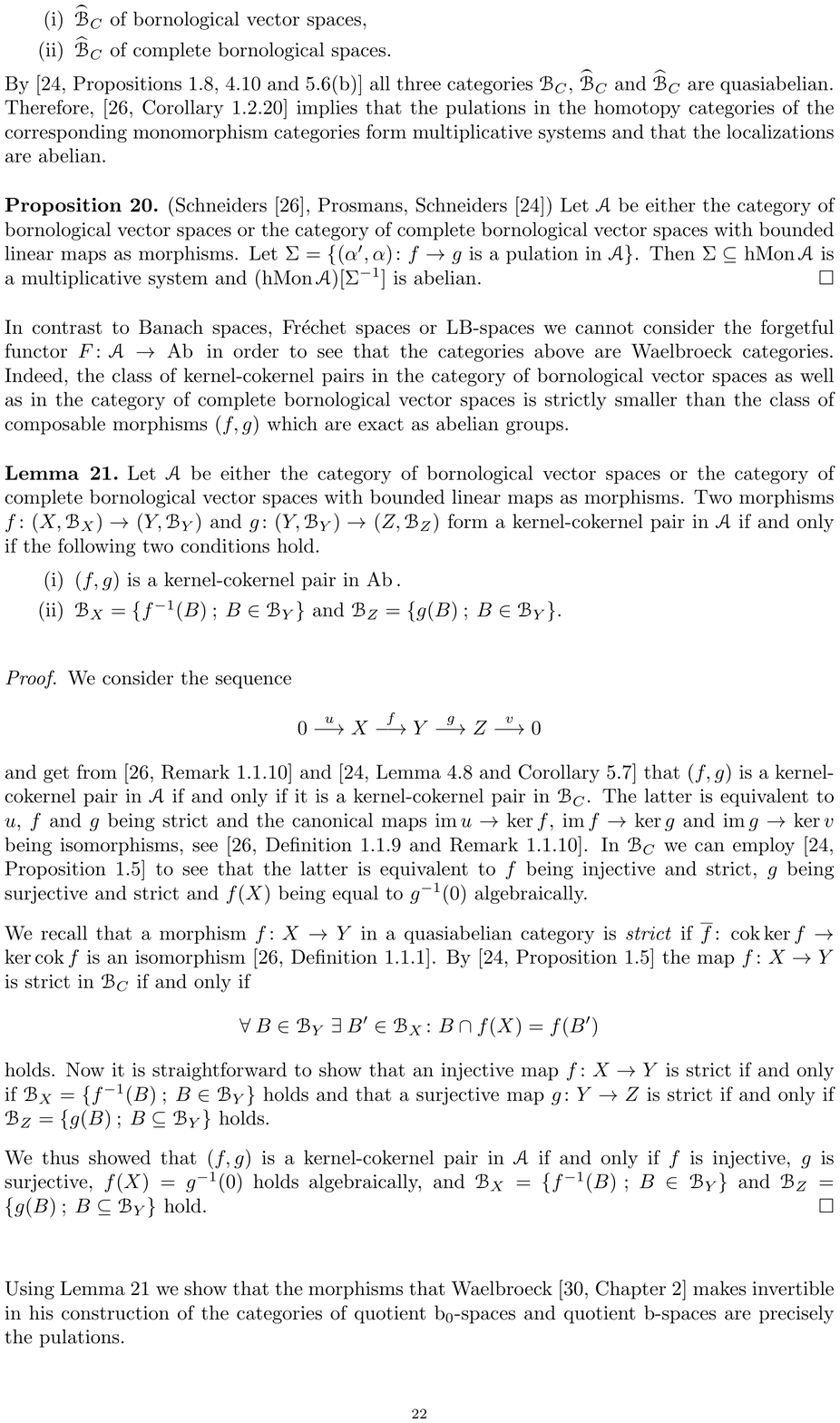, width=21cm}}
\end{picture}
}
\newpage
\parbox{1000pt}{
\begin{picture}(0,0)(45,629)
\put(89,-2){\psfig{figure=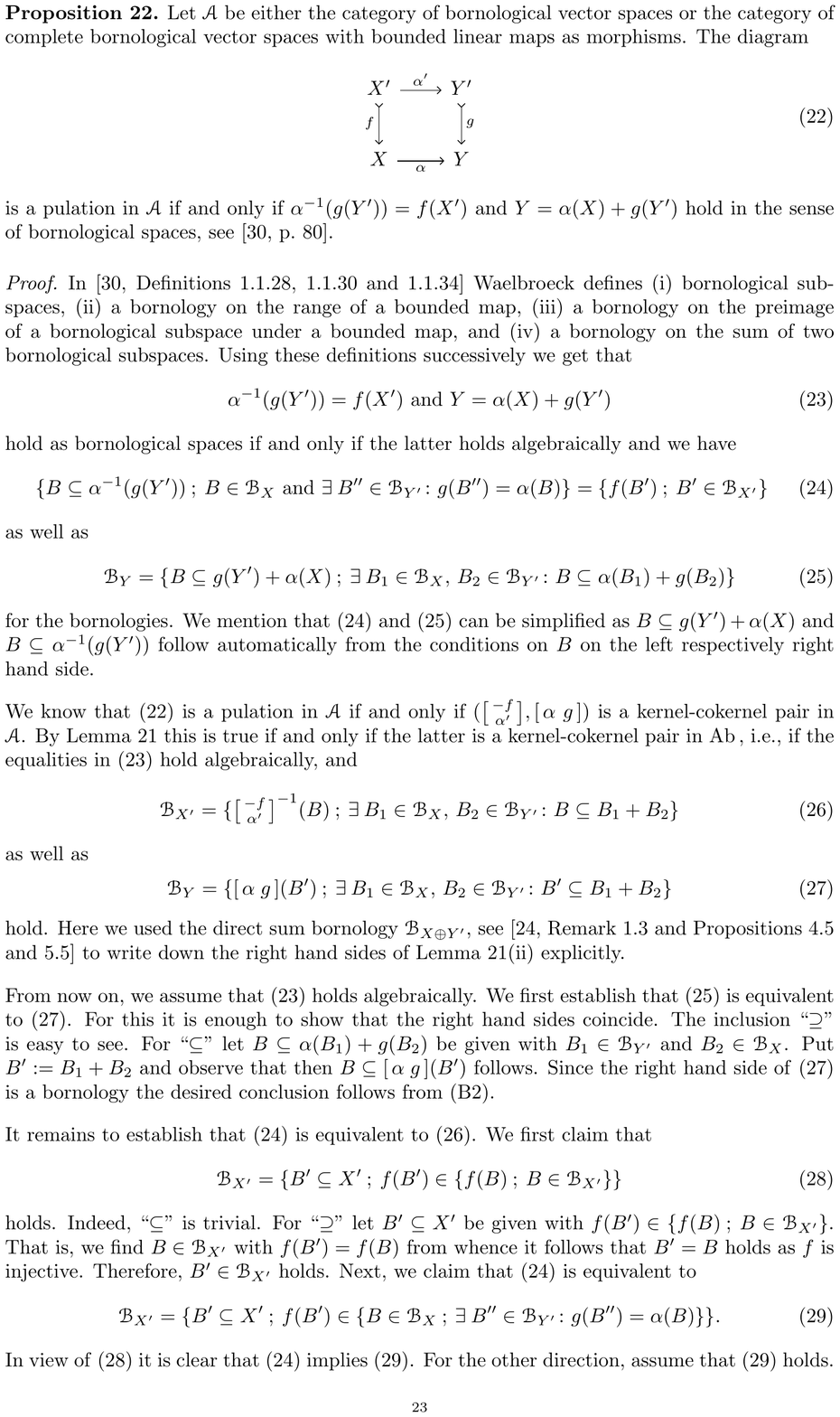, width=21cm}}
\end{picture}
}
\newpage
\parbox{1000pt}{
\begin{picture}(0,0)(45,629)
\put(89,-2){\psfig{figure=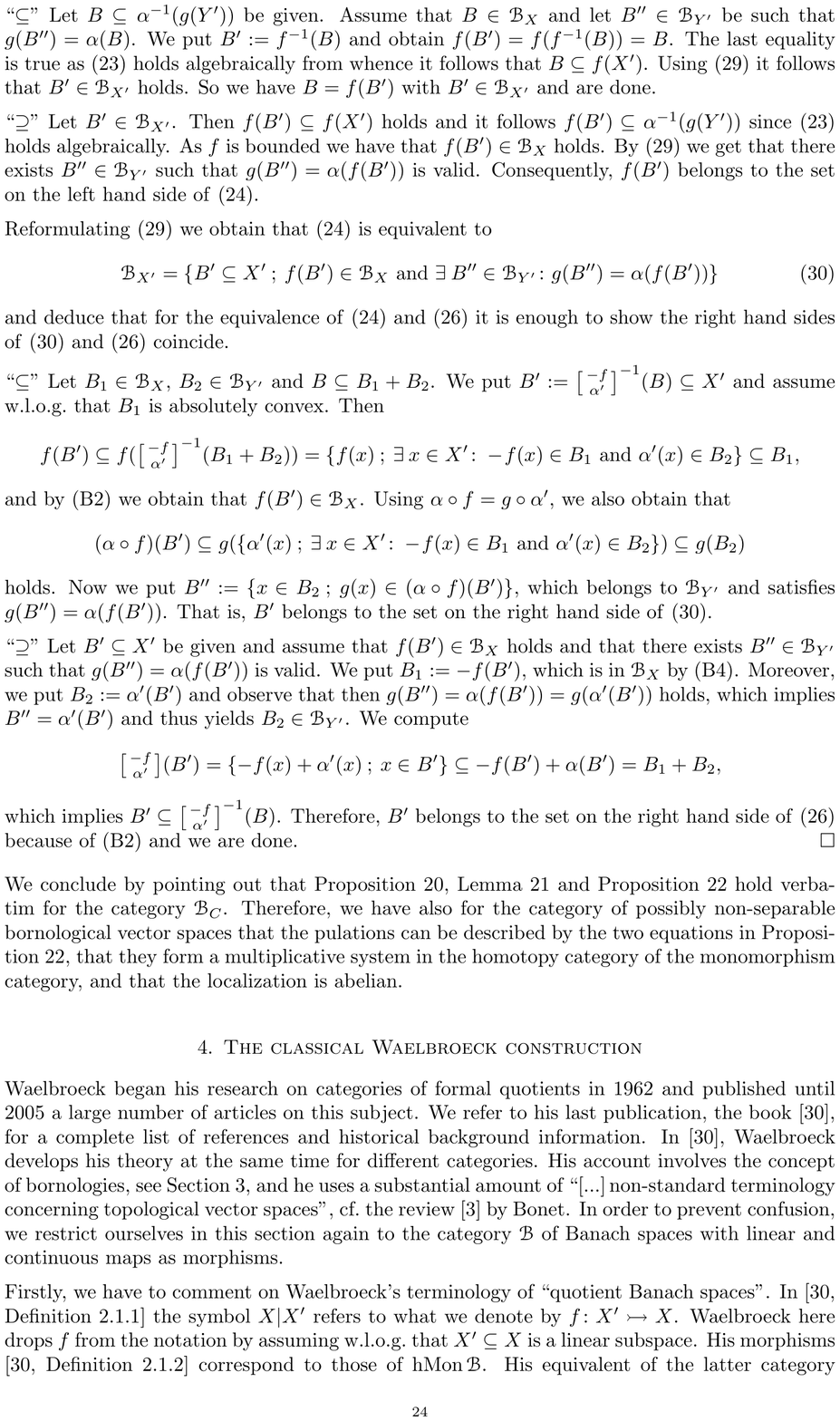, width=21cm}}
\end{picture}
}
\newpage
\parbox{1000pt}{
\begin{picture}(0,0)(45,629)
\put(89,-2){\psfig{figure=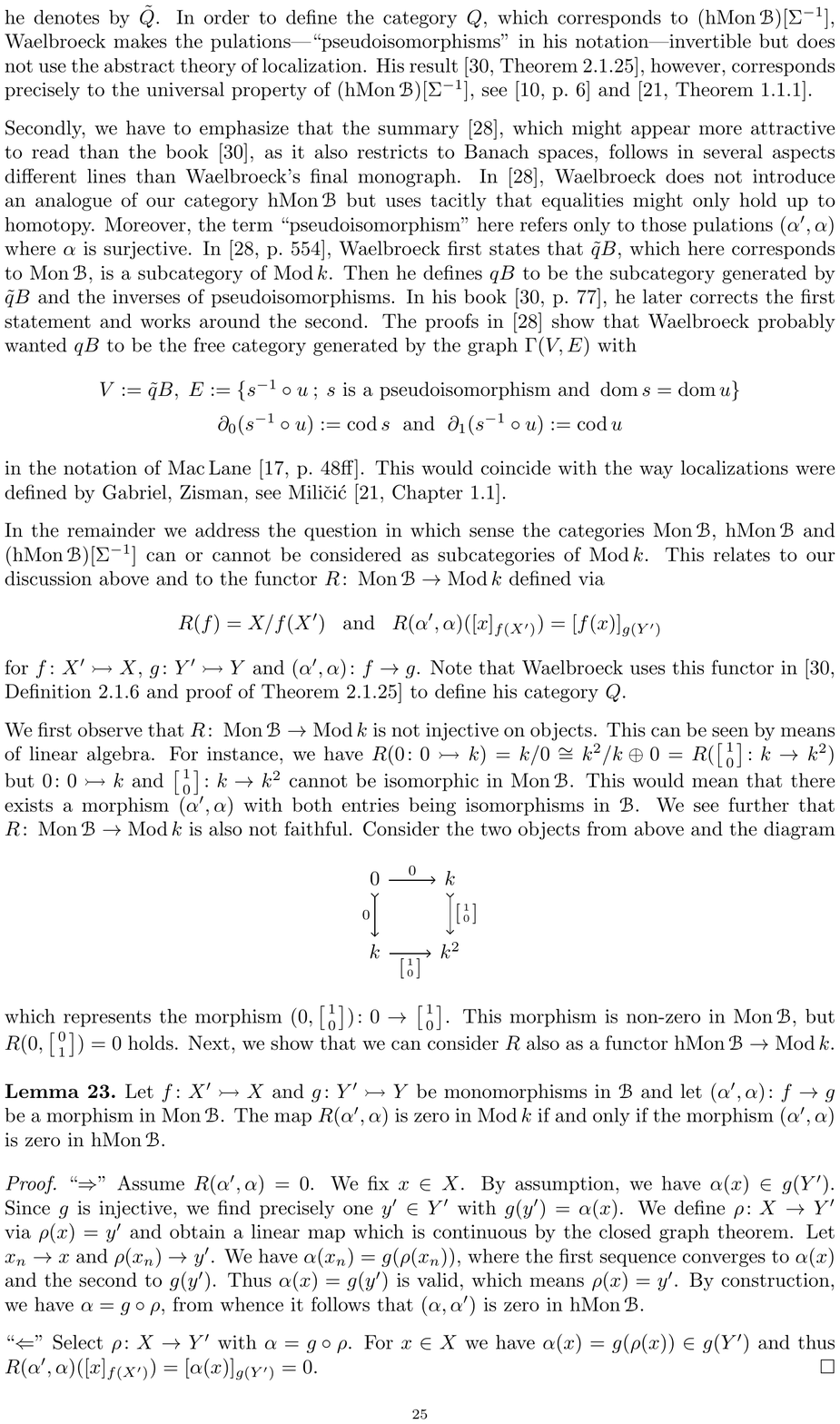, width=21cm}}
\end{picture}
}
\newpage
\parbox{1000pt}{
\begin{picture}(0,0)(45,629)
\put(89,-2){\psfig{figure=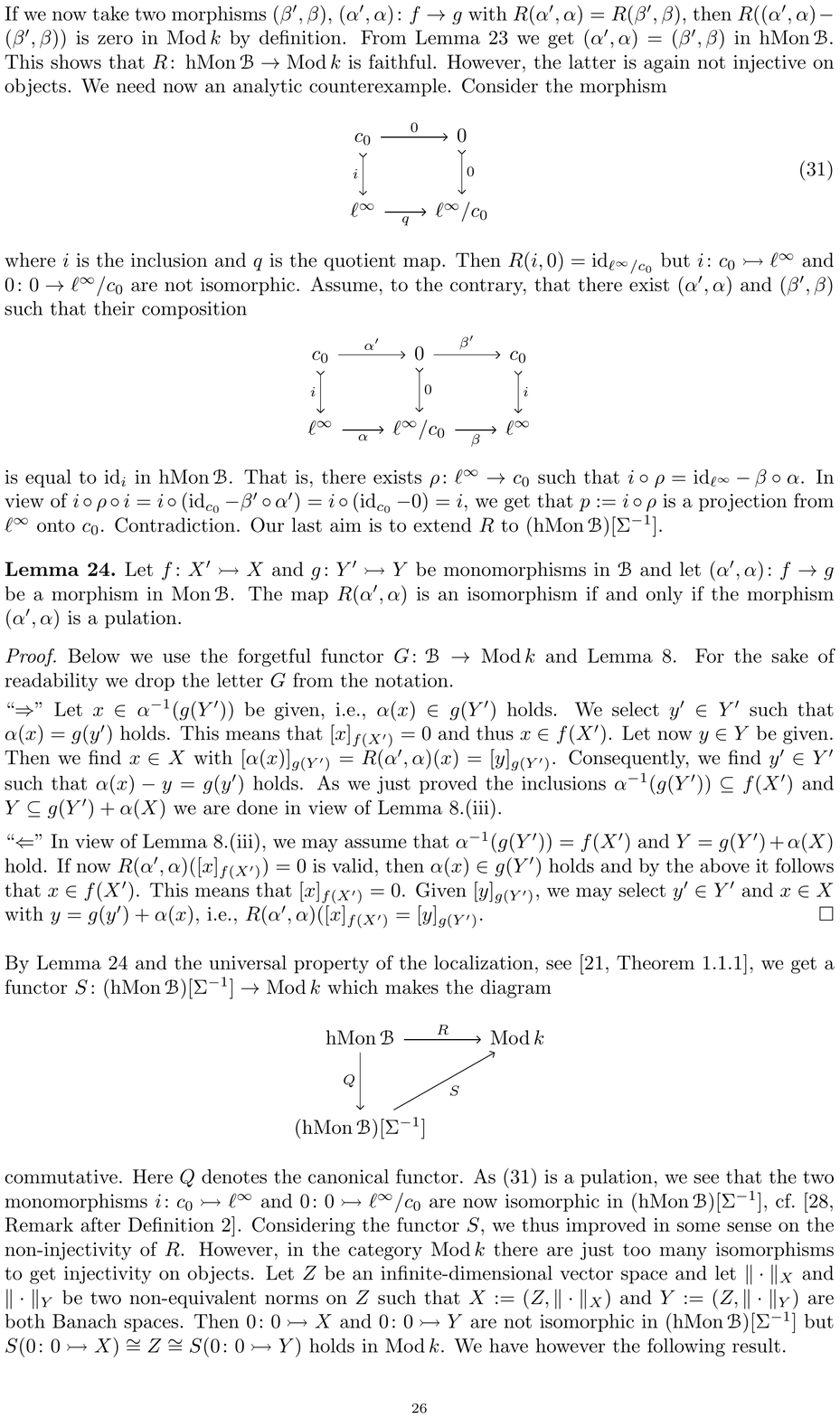, width=21cm}}
\end{picture}
}
\newpage
\parbox{1000pt}{
\begin{picture}(0,0)(45,629)
\put(89,-2){\psfig{figure=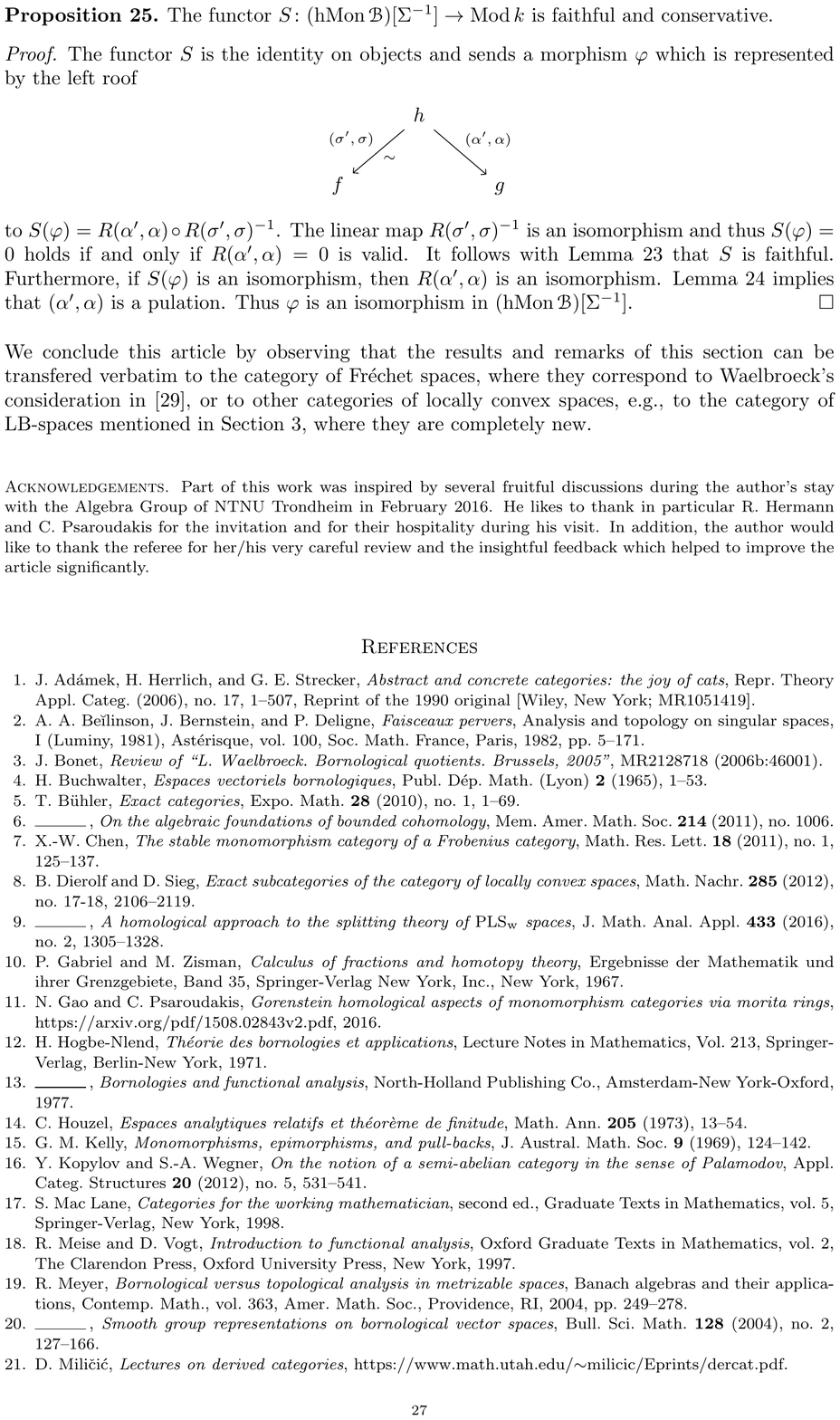, width=21cm}}
\end{picture}
}
\newpage
\parbox{1000pt}{
\begin{picture}(0,0)(45,629)
\put(89,-2){\psfig{figure=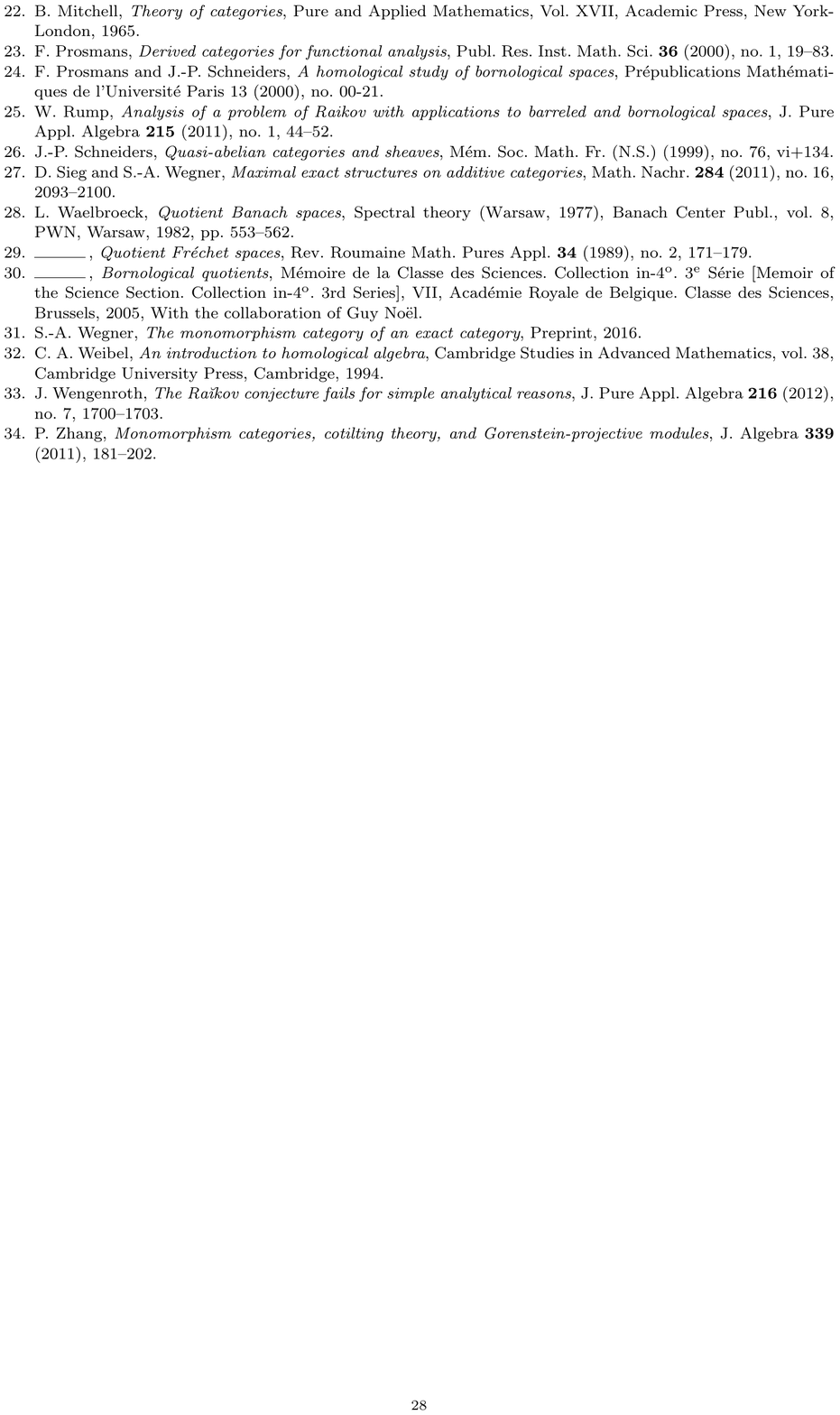, width=21cm}}
\end{picture}
}

\end{document}